\documentclass{amsart}
\usepackage{amsmath,amssymb}
\newtheorem{Theorem}{Theorem}[section]

\newtheorem{Lemma}[Theorem]{Lemma}
\newtheorem{Proposition}[Theorem]{Proposition}
\theoremstyle{definition}
\newtheorem{Definition}[Theorem]{Definition}
\theoremstyle{remark}
\newtheorem{rem}[Theorem]{Remark}
\theoremstyle{remark}

\numberwithin{equation}{section}
\def\K{{ \! \rm \ I\!K}}

\newcommand{\R}{\mathbb R}
\newcommand{\C}{\mathbb C}
\newcommand{\N}{\mathbb N}
\newcommand{\Z}{\mathbb Z}
\newcommand{\eps}{\varepsilon}

\newcommand{\Con}{\mathcal{C}}
\newcommand{\F}{\mathcal{F}}
\newcommand{\pol}{\mathfrak{Pol}}
\newcommand{\normal}{\mathcal{N}}
\newcommand{\res}{\hbox{\rm res}}
\newcommand{\tr}{\hbox{\rm tr}}

\begin{document}

\title{On Diff(M)-pseudo-differential operators and the geometry of non linear grassmannians}%
\author{Jean-Pierre Magnot}%
\address{Lyc'ee Jeanne d'Arc, Avenue de Grande Bretagne, F-63000 Clermont-Ferrand}%
\email{jean-pierr.magnot@ac-clermont.fr}%

\begin{abstract}
We consider two principal bundles of embeddings with total space $Emb(M,N),$ with structure groups $Diff(M)$
and $Diff_+(M),$ where $Diff_+(M)$ is the groups of orientation preserving diffeomorphisms.
The aim of this paper is to describe the structure group of the tangent bundle of the two base manifolds: 
$$ B(M,N) = Emb(M,N)/Diff(M) \hbox{ and } B_+(M,N)= Emb(M,N)/Diff_+(M).$$
From the various properties described, an adequate group seems to be a group of Fourier integral operators, which is carefully studied. 
This is the main goal of this paper to analyze this group, which is a central extension of a group of diffeomorphisms by a group of pseudo-differential operators which is slightly different from the one developped in \cite{OMYK4}.  We show that these groups are regular, and develop the necessary properties for applications to the geometry of $ B(M,N) .$
A case of particular interest is $M=S^1,$ where connected components of $B_+(S^1,N)$ are deeply linked with homotopy classes of oriented knots. In this example, the structure group of  the tangent space $TB_+(S^1,N)$ is a subgroup of some group $GL_{res},$ following the classical notations of \cite{PS}. These constructions  suggest 
some 	approaches in the spirit of \cite{Ma2006} that could lead to knot invariants through a theory of Chern-Weil forms.  
\end{abstract}

\maketitle



\tableofcontents

\section*{Introduction}

Given $M$ and $N$ two Riemannian manifolds without boundary, with $M$ compact, 
the space of smooth embeddings $Emb(M,N)$ is currently known as a principal bundle with structure group $Diff(M),$ where $Diff(M)$ naturally acts by composition of maps. 
The base $$B(M,N)=Emb(M,N)/Diff(M)$$ is known as a Fr\'echet manifold, and there exists some local trivializations of this bundle. 
We focuse here on the base manifold, which seems to carry a richer structure than $Emb(M,N)$ itself. 

This paper gives the detailed description of the structure group of the tangent bundle of connected components of $TB(M,N).$ This
structure group can be slightly different when changing of connected component of $B(M,N).$ 
It is viewed as an extension of the group of automorphisms $Aut(E)$ of a vector bundle $E$ by some group of pseudo-differential operators.  
We show that this group is a regular Lie group (in the sense that it carries an exponential map),
and that it is also a group of Fourier integral operators, which explains the notations $FIO_{Diff}$ and $FCl_{Diff}$ (``$Cl$'' for ``classical''). 
All these groups are constructed along a short exact sequence of the type 
$$`` 0 \rightarrow PDO \rightarrow FIO \rightarrow Diff \rightarrow 0'',$$
where $PDO$ is a group of pseudo-differential operators, FIO is a group of Fourier integral operators, and $Diff$ is a group of diffeomorphisms;
this sequence plays a central role in the proofs.
The theorems described are general enough to be applied to many groups of diffeomorphisms: volume preserving diffeomorphisms, symplectic diffeomorphisms, hamiltonian diffeomorphisms, and to groups of pseudo-differential operators: classical or non-classical, bounded or unbounded, compact and so on, but we concentrate our efforts on $Diff(M)$ and $Diff_+(M),$ the group of orientation preserving diffeomorphisms. The constructions are made for operators acting on smooth sections of trivial or non trivial bundles. For a non trivial bundle $E,$
 the group of automorphisms of the bundle plays a central role in the description, because easy arguments suggest that there is no adequate embedding of the group of diffeomorphisms of the base manifold into the group of automorphisms of the bundle
Specializing to $M=S^1,$
given a (real) vector bundle $E$ over $S^1,$ the groups  $FIO_{Diff}(S^1,E) $ and in particular  $FCl^{0,*}_{Diff_+}(S^1,E)$ is of particular interest, where $FCl^{0,*}_{Diff_+}(S^1,E)$
is defined through the short exact sequence:
$$ 0 \rightarrow Cl^{0,*}(S^1,E) \rightarrow FCl^{0,*}_{Diff_+}(S^1,E)\rightarrow Diff_+(S^1) \rightarrow 0,$$
where $Cl^{0,*}(S^1,E)$ is the group of bounded classical pseudo-differential operators and $Diff_+(S^1)$ is the group of orientation-preserving diffeomorphisms. We have to notice that the necessary Fourier analysis on these operators naturally takes place in the complexification $E_\C$ of the vector bundle $E$, and that $E_\C$ as a complex vector bundle is trivial, but the real vector bundle $E$ can be non trivial. 
Given any Riemannian connection on the bundle $E$, if $\epsilon$ is the sign of this connection (and this is a bounded pseudo-differential operators acting on smooth sections of $E$), it appears that $[ FCl^{0,*}_{Diff_+}(S^1,E), \epsilon] $ is a set of smoothing operators. 
Thus, it is a subgroup of the group $$Gl_{res} = \left\{ u \in Gl(L^2(S^1,E)) | [\epsilon, u]  \hbox{ is Hilbert-Schmidt } \right\}.$$
Even if the inclusion is not a bounded inclusion, this result extends the results given in \cite{PS} on the group $Diff_+(S^1)$
(which inclusion map into $Gl_{res}$ is not bounded too) and in \cite{Ma2006-2} for the group $Cl^{0,*}(S^1,E).$ 
We get a non-trivial cocycle on the Lie algebra of $FCl^{0,*}_{Diff_+}(S^1,E)$
by the Schwinger cocycle, extending results obtained in \cite{Ma2006-2,Ma2008} for a trivial complex bundle. 
 
Coming back to $Emb(M,N),$ 
one could suggest that $Aut(E)$ is sufficient as a structure group, but we refer the reader to earlier works such as \cite{F2,CDMP,Ma2006}
to see how pseudo-differential operators can arise from Levi-Civita connections of Sobolev metrics 
when the adequate structure group for the $L^2$ metric is a group of multiplication operators. 
Moreover, especially for $M=S^1,$ taking the quotient $$B_+(S^1,N) = Emb(S^1,N)/Diff_+(S^1),$$
we show that there is a sign operator $\epsilon(D)$, which is a pseudo-differential operator of order 0, and coming intrinsically from the geometry of  $Emb(S^1,N),$ 
such that the recognized structure group of $TB_+(S^1,N)$ is $FCl^{0,*}_{Diff_+}(S^1,E)\subset Gl_{res}.$  

As a motivation for our description, which was inspired by exchanges with Claude Roger, we give some Chern-Weil type forms 
along the lines of \cite{Ma2006}, that is, differential forms of the type 
$tr(\Omega^k).$ In the proposed application, $``tr''$ is the Kontsevich-Vishik 
trace of odd class pseudo-differential operators \cite{KV1,KV2}. When $M$ is odd-dimensional, 
these Chern-Weil forms are closed and do not depend on the choice of the 
connection. They can be non exact on non trivial principal bundles.  So that, 
over each 
connected component $C$ of $B(M,N),$
we have a principal bundle $Emb(M,N)|_C.$ Chern-Weil forms $tr(\Omega^k)$ 
describe obstructions to the triviality of the principal bundle $Emb(M,N)|_C.$
From another point of view, when $M=S^1,$ $C$ describes the smooth homotopy class of a knot. For an arbitrary manifold $M,$ some authors talk about ``generalized knots''. These Chern-Weil forms, among others, can describe knot invariant. This is a work actually in preparation.

\section{Preliminaries on algebras and groups of operators}

Now, we fix $M$ the source manifold, which is assumed to be Riemannian, 
compact, connected and without boundary, and the target manifold 
which is only assumed Riemannian. We note by $Vect(M)$ the space of vector fields on $TM.$ 
Recall that the Lie algebra of the group of diffeomorphisms is $Vect(M),$ which is a Lie-subalgebra of the (Lie-)algebra of 
differential operators, which is itself a subalgebra of the algebra of 
classical pseudo-differential operators.

\subsection{Differential and pseudodifferential operators on a manifold $M$}

 \begin{Definition}
Let $DO(M)$ be the graded algebra of operators, acting on $C^\infty(M,\R)$, generated by:

$\bullet$ the multiplication operators: for $f \in C^\infty(M,\R),$ we define the multiplication operator 
$$M_f : g \in C^\infty(M,\R) \mapsto f.g \hbox{ (by pointwise multiplication)}$$

$\bullet$ the vector fields on $M$: for a vector field $X \in Vect(M)$, we define the 
differentiation operator
$$D_X : g \in C^\infty(M,\R) \mapsto D_X g \hbox{ (by differentiation, pointwise)}$$
Multiplication operators are operators of order $0$, vector fields are operators of order 1. 
For $k \geq 0,$ we note by $DO^k(M)$ the differential operators of order $\leq k.$
\end{Definition}
        
Differential operators are \textbf{local}, which means that
$$ \forall A \in DO(M), \forall f \in C^\infty_c(M,\mathbb{R}), supp(A(f)) \subset supp(f).$$
The inclusion $Vect(M)\subset DO(M)$ is an inclusion of Lie algebras.
The algebra $DO(M),$ graded by the order, is a subalgebra of the algebra of classical pseudo-differential operators
$Cl(M),$ which is an algebra that contains the square root of the Laplacian, and its inverse. 
This algebra contains trace-class operators on $L^2(M,\R).$ 
An exposition  of basic facts on pseudo-differential operators defined 
on a vector bundle $E \rightarrow M$ can
be found in \cite{Gil} for definition of pseudo-differential operators 
and of their order, 
(local) definition of symbols and spectral properties. 
We assume known the definition of the algebra of pseudo-differential operators
$PDO(M,E)$,  
classical pseudo-differential operators $Cl(M,E)$. When the vector bundle $E$ is assumed trivial, i.e. $E = M\times V$ or $E = M \times \K^p$ with $\K = \R$ or $\C,$ we use the notation $Cl(M,V)$ or $Cl(M,\K^p)$ instead of $Cl(M,E).$ These operators are \textbf{pseudolocal}, which means that 
$$ \forall A \in PDO(M,E), \forall f \in L^2(M,E), \hbox{ if } f \hbox{ is smooth on } K, \hbox{ then } A(f) \hbox{ is smooth on } K.$$ 
\begin{Definition}
A pseudo-differential operator $A$ is \textbf{log-polyhomogeneous}
if and only if its formal symbol reads (locally) as
$$\sigma(A)(x,\xi)\sim_{|\xi| \rightarrow +\infty} \sum_{j=0}^o \sum_{k = -\infty}^{o'}
\sigma_{j,k}(x,\xi)(log(|\xi|))^j,$$
where $\sigma_{j,k}$ is a positively $k-$homogeneous symbol.
\end{Definition}
The set of log-polyhomogenous pseudo-differential operators 
is an algebra.  

\vskip 12pt
A global symbolic calculus has been defined independently 
by two authors in \cite{BK}, \cite{Wid}, where we can see how the 
geometry of the base manifold $M$ furnishes an obstruction to generalize 
local formulas of composition ans inversion of symbols. We do not recall these
formulas here because they are not involved in our computations. 
More interesting for this article is to precise when the local formulas 
of composition of formal symbols extend globally on the base manifold. 

We assume that $M$ is equipped with charts  such that the changes 
of coordinates are translations and that the vector bundle $E \rightarrow M$ 
is trivial.  
This is in particular true when  
$M = S^1 = {\R \over 2\pi\Z},$ or when $M = T^n = \prod_{i=1}^n S^1.$ 
In the case of $S^1,$ we use the smooth atlas 
$\mathcal{ATL}$ of $S^1$ defined as follows:
\begin{eqnarray*}
\mathcal{ATL} & = & \{\varphi_0,\varphi_1\} ; \\
\varphi_n & : & x \in ]0 ; 2\pi[ \mapsto e^{i(x + n\pi)} \subset S^1 
\hbox{ for } n\in \{0;1\} \end{eqnarray*}     
Associated to this atlas, we fix a smooth partition of the unit $\{s_0;s_1\}$.
An operator $A : C^\infty(S^1, \C) \rightarrow C^\infty(S^1,\C)$ can be described in 
terms of 4 operators
$$ A_{m,n} : f \mapsto s_m \circ A \circ s_n \hbox{ for } (m,n) \in \{0,1\}.$$

Such a formula is a straightforward application of a localization formula in the case of an atlas $\{ \varphi_i \}_{i \in I}$ of a manifold $M$ with associated family of partitions of the unit $\{s_i \}_{i \in I},$ see e.g \cite{Gil} for details.

\vskip 12pt
\noindent
\textbf{Notations.} 
We note by  $ PDO (M, \C) $
(resp.  $ PDO^o (M, \C)
$, resp. $Cl(M,\C)$) the space of
pseudo-differential operators (resp.
pseudo-differential operators of order o, resp. classical pseudo-differential operators) acting on smooth
sections of $E$, and by $Cl^o(M,\C)= PDO^o(S^1,\C) \cap Cl(S^1,\C)$ the space of classical 
pseudo-differential operators of order $o$.
\vskip 10pt

If we set
$$ PDO^{-\infty}(M,\C) = \bigcap_{o \in \Z} PDO^o(M,\C),$$
we notice that it is a two-sided ideal of $PDO(M,\C)$, and we define the quotient algebra
$$\mathcal{F}PDO(M, \C) = PDO(M,\C) / PDO^{-\infty}(M,\C),$$ 
$$\F Cl(M,\C) = Cl(M,\C) /  PDO^{-\infty}(M,\C),$$
$$ \quad \F Cl^o(M,\C) = Cl^o(M,\C)  / PDO^{-\infty}(M,\C),$$
called the algebras of formal pseudo-differential operators. 
$\mathcal{F}PDO(M,\C)$ is isomorphic to the set of formal symbols \cite{BK}, 
and the identification is a morphism of $\C$-algebras, for the multiplication on formal symbols
defined before (see e.g. \cite{Gil}). At the level of kernels of operators, a smoothing operator has a kernel $K_\infty \in C^\infty(M\times M, \C)$, where as the kernel of a pseudo-differential operator is in general smooth only on the off-diagonal region $(M\times M)-\Delta(M),$ where $\Delta(M)$ denotes here, very exceptionnally in this paper, the diagonal set (and not a Laplacian operator). We finish by mentonning that the last property is equivalent to pseudo-locality.  

\subsection{Fourier integral operators}
With the notations that we have set before, a scalar
Fourier-integral operator of order $o$ 
is an operator $$ A : C^\infty(M, \C) \rightarrow C^\infty(M,\C)$$
such that, $\forall (i,j) \in I^2,$
\begin{eqnarray} \label{localization} A_{k,j}(f) &  = & \int_{supp(s_j)} e^{-i\phi(x,\xi)}\sigma_{k,j}(x,\xi) \hat{(s_j.  
f)} (\xi) d\xi \end{eqnarray}
where $\sigma_{k,j} \in C^\infty( supp(s_j) \times \R, \C)$ satisfies
$$\forall (\alpha, \beta) \in \N^2, \quad |D^\alpha_x D^\beta_\xi \sigma_{k,j}(x,\xi)
|\leq C_{\alpha,\beta}(1 + |\xi|)^{o-\beta},$$
and where, on any domain $U$ of a chart on $M$, 
$$\phi(x,\xi) : T^*U-U\approx U\times \R^{dim M}_\{0\}\rightarrow \R$$
is a smooth map, positively homogeneous of degree 1 fiberwise and such that $$\det\left(\frac{\partial^2\phi }{\partial x \partial_\xi} \right) \neq 0.$$
Such a map is called \textbf{phase} function.
(In these formulas, the maps are read on local 
charts but we preferred to only mention this aspect and not 
to give heavier formulas and notations) An operator $A$ is  {pseudo-differential} operator
if the operators $A_{k,l}$ in formula \ref{localization} can be written 
as Fourier integral operators 
with $\varphi(x,\xi) = x.\xi.$ 
Notice that, in order to define an operator $A,$ the choice of 
$\varphi$ and $\sigma_{k,l}$ is not a priori unique for 
general Fourier integral operators.
Let $E = S^1 \times \C^k$ be a trivial smooth vector bundle over $S^1$.
An operator acting on $C^\infty(M,\C^n)$ is Fourier integral operator 
(resp. a pseudo-differential operator) 
if it can be viewed as a $(n \times n)$-matrix of Fourier integral operators with same phase function 
(resp. scalar pseudo-differential operators).

We define also the algebra of formal operators, which is the quotient space
$$\F FIO = FIO / PDO^{-\infty},$$ which is possible because $PDO^{-\infty}$
is a closed two-sided ideal. When we consider 
classical Fourier integral operators, noted $FCl,$ that is  
operators with classical symbols, we add
to this topology the topology on formal symbols \cite{ARS1,ARS2}
which is an ILH topology (see e.g. \cite{Om} for state of the art). 
We want to quote that if the symbols $\sigma_{m,n}$ 
are symbols of order $0,$ then we 
get Fourier integral operators that are $L^2-$bounded. We note
this set $FIO^0.$ This set is a subset of $FIO,$ and we have 
$$Cl^0 \subset PDO^0 \subset FIO^0 \subset FIO.$$

The techniques used for pseudo-differential operators are also used 
on Fourier integral operators, especially Kernel analysis.
Let us consider a local coordinate operator $A_{m,n}$
then, using the notation of of the formula \ref{localization}, the operator
$A_{m,n}$is described by a kernel 
$$K_{m,n}(x,y) = \int_\xi e^{-i\left(\phi\left((x),\xi)\right)-y.\xi\right)} \sigma_{m,n}(x,\xi) d\xi.$$
From this approach one derives the composition and inversion formulas that 
will not be used in this paper, see e.g. \cite{Horm},
but in the sequel we shall use the slightly restricted class
of operators studied in \cite{OMY1,OMY2,OMYK3,OMYK4,OMYK5,OMYK6,OMYK7,OMYK8} 
and also in \cite{ARS1,ARS2} for formal operators. 

\subsection{Topological structures and regular Lie groups of operators}

The topological structures can be derived both from 
symbols and from kernels, as we have quoted before
but principally because there is the exact sequence described below with slice.
At the level of units of these sets, i.e. of groups of invertible operators, 
the existence of the slice is also crucial.  
In the papers \cite{ARS1,ARS2,OMY1,OMY2,OMYK3,OMYK4,OMYK5,OMYK6,OMYK7,OMYK8,RS}, 
the group of invertible Fourier integral operators receives first 
a structure of topological group, 
with in addition a differentiable structure, e.g. a Fr\"olicher structure,
which recognized as a structure of generalized Lie group, see e.g. \cite{Om}. 

We have
to say that, with the actual state of knowledge, using \cite{KM}, we can give  
a manifold structure (in the convenient setting described by Kriegl and Michor or in the category of Fr\"olicher spaces following \cite{Ma2013}) 
to the corresponding Lie groups.
Let us recall the statement

\begin{Theorem}\label{slice} \cite{KM}
Let $G,H,K$ be convenient Lie groups or Fr\"olicher Lie groups such that there is a short exact sequence of Lie groups
$$ 0 \rightarrow H \rightarrow G \rightarrow K \rightarrow 0$$
such that there is a local slice $ K \rightarrow G$. Then 
$$ G \hbox{ regular } \Leftrightarrow H \hbox{ and } K \hbox{ regular.}$$
\end{Theorem}

\begin{rem}
In \cite{ARS1,ARS2,OMY1,OMY2,OMYK3,OMYK4,OMYK5,OMYK6,OMYK7,OMYK8,RS}, the group $K$ considered is the group of 1-positively homogeneous symplectomorphisms $Diff_\omega(T^*M - M)$ where $\omega$ is the canonical symplectic form on the cotangent bundle. The local section considered enables to build up the phase function of a Fourier integral operator from such a symplectic diffeomorphism inside a neighborhood of $Id_M.$ There is a priori no reason to restrict the constructions to classical pseudo-differential operators of order 0, and have groups to Fourier integral operators with symbols  in wider classes. This remark appears important to us because the authors cited before restricted themselves to classical symbols.
\end{rem}

\subsection{$PDO(M,E),$ $Aut(E)$ and $Diff(M).$}
We get now to another group:
\begin{Theorem} \label{FIOdiff}
Let $H$ be a regular Lie group of pseudo-differential operators acting on smooth sections of  a trivial bundle $E \sim V \times M \rightarrow M.$
The group $Diff(M)$ acts smoothly on  $C^\infty(M,V),$ and is assumed to act smoothly on $H$ by adjoint action.
If $H$ is stable under the $Diff(M)-$adjoint action, then there exists a corresponding regular Lie group $G$ of Fourier integral operators through the exact sequence:
$$ 0 \rightarrow H \rightarrow G \rightarrow Diff(M) \rightarrow 0.$$ If $H$ is a Fr\"olicher Lie group, then $G$ is a Fr\"olicher Lie group. If $H$ is a Fr\'echet Lie group, then $G$ is a Fr\'echet Lie group.
\end{Theorem}

\begin{rem}
The pseudo-differential operators can be classical, log-polyhomogeneous, or anything else. Applying the formulas of ``changes of coordinates'' (which can be understood as adjoint actions of diffeomorphisms) of e.g. \cite{Gil}, one easily gets the result.
\end{rem}

\vskip 12pt
\noindent 
\textbf{Proof of Theorem \ref{FIOdiff}.}
Let us first notice that the action $$ (f,g) \in C^\infty(M,V) \times Diff(M) \mapsto f \circ g \in C^\infty(M,V)$$
can be read as, first a linear operator $T_g$ with kernel $$K(x,y) = \delta(g(x),y) \quad \hbox{(Dirac} \delta-\hbox{function)}$$
or equivalently, on an adequate system of trivializations \cite{Gil},
$$ T_g(f)(x) = \int e^{ig(x).\xi} \hat{f}(\xi) d\xi.$$
This operator is not a pseudo-differential operator because it is not pseudolocal (unless $g = Id_M$), but since $$det \left(\partial_x\partial_\xi (g(x).\xi) \right)= det (D_xg),$$ we get that $T_g$ is a Fourier-integral operator. Notice that another way to see it is the expression of its kernel. 

Now, given $(A, g) \in H \times Diff(M),$ we define $$A_g = T_g \circ A.$$ We get here a set $G$ of operators which is set-theorically isomorphic to $H \times Diff(M).$
Since $H$ is invariant under the adjoint action of the group $Diff(M),$ G is a group, and from the beginning of this proof, we get that $G$ is a group, and that there is the short exact sequence announced:
 $$ 0 \rightarrow H \rightarrow G \rightarrow Diff(M) \rightarrow 0,$$
with a global slice $$g \in Diff(M) \mapsto T_g \in G.$$
Since the adjoint action of Diff(M) is assumed smooth on $H$, we can endow $G$ with the product Fr\"olicher structure to get a regular Fr\"olicher Lie group. Since $Diff(M)$ is a Fr\'echet Lie group, if $H$ is a Fr\'echet Lie group, then $G$ is a Fr\'echet Lie group. \qed  
 
\begin{rem}
Some restricted classed of such operators are already considered in the literature under the name of $G-$pseudo-differential operators, see e.g.\cite{SS}, but the groups considered are discrete (amenable) groups of diffeomorphisms. 
\end{rem}
 \begin{Definition}
Let $M$ be a compact manifold and $E$ be a (finite rank) trivial vector bundle over $M.$
We define $$ FIO_{Diff}(M,E) = \left\{ A \in FIO(M,E) | \phi_A(x,\xi) = g(x).\xi; g \in Diff(M) \right\}.$$
\end{Definition}
The set of invertible operators $FIO_{Diff}^*(M,E)$ is obviously a group, that decomposes as 
$$ 0 \rightarrow PDO^*(M,E) \rightarrow FIO_{Diff}^*(M,E) \rightarrow Diff(M) \rightarrow 0$$
with global smooth section $$ g \in Diff(M) \mapsto (f \in C^\infty(S^1,E) \mapsto f \circ g ).$$
Hence, Theorem \ref{FIOdiff} applies trivially to the following context:

\begin{Proposition}
Let $FCl^{0,*}_{Diff}(M,E) $ be the set of operators $ A \in FIO_{Diff}^*(M,E)$ such that $A$ has a 0-order classical symbol. Then we get the exact sequence:
$$ 0 \rightarrow Cl^{0,*}(M,E) \rightarrow FCl_{Diff}^{0;*}(M,E) \rightarrow Diff(M) \rightarrow 0$$
and $FCl^{0,*}_{Diff}(M,E) $ is a regular Fr\"olicher Lie group, with Lie algebra isomorphic, as a vector space, to $Cl^0(M,E) \oplus Vect(M).$
\end{Proposition}
 
Notice that the triviality of the vector bundle $E$ is here essential to make a $Diff(M)-$action on smooth section of $C^\infty(M,E).$
Let us assume now that $E$ is not trivial. 
At the infinitesimal level, trying to extend straightway, one gets a first condition for the extension.

\begin{Lemma} (see e.g. \cite{AMR2004})
Let us fix a $0-$curvature connection $\nabla$ on $M.$ Then $X \in Vect(M) \mapsto \nabla_X \in DO^1(M,E)$ is a one-to-one Lie algebra morphism.
\end{Lemma}

We remark that the analogy with the setting of trivial bundles $E$ stops here since the group $Diff(M)$ 
cannot be recovered in this group of operators. For example, when $M=S^1$, if $E$ is non trivial, 
the (infinitesimally) flat connection ensures that the holonomy group $\mathcal{H}$ is discrete, but it cannot be trivial since 
the vector bundle $E$ is not. 
On non a trivial bundle $E,$ let us consider the group of bundle automorphism $Aut(E).$ 
The gauge group $DO^0(M,E)$ is naturally embedded in $Aut(E)$ and the bundle projection $$E \rightarrow M$$ 
induces a group projection $$\pi:Aut(E) \rightarrow Diff(M).$$
Therefor we get a short exact sequence 
$$ 0 \rightarrow DO^{0,*}(M,E) \rightarrow Aut(E) \rightarrow Diff(M) \rightarrow 0.$$
Following \cite{ACMM1989} there exists a local slice $U \subset Diff(M) \rightarrow Aut(E),$ where $U$ is a $C^0-$open neighborhood on $Id_M,$ 
which shows that $Aut(E)$ is a regular Fr\'echet Lie group.
Therefore, the smallest group spanned by $PDO^*(M,E)$ and $Aut(E)$ is such that:
\begin{itemize}
\item the projection $E \rightarrow M$ induces a map $Aut(E) \rightarrow Diff(M)$ with kernel $DO^0(M,E) = Aut(E) \cap PDO(M,E)$
\item $Ad_{Aut(E)}(PDO(M,E)) = PDO(M,E)$  
\end{itemize} 
therefore we can consider the space of operators on $C^\infty(M,E)$
$$ FIO_{Diff}^*(M,E) = Aut(E) \circ PDO^*(M,E).$$

\begin{Lemma} \label{phase}
The map $$ (B,A) \in Aut(E) \times PDO^*(M,E) \mapsto \pi(B) \in Diff(M)$$ induces a ``phase map'' $$\tilde{\pi}:FIO_{Diff}^*(M,E) \rightarrow Diff(M).$$
\end{Lemma}

\noindent
\textbf{Proof.}
Let $((B,A),(B',A')) \in Aut(E) \times PDO^*(M,E).$
\begin{eqnarray*} B\circ A = B' \circ A' &\Leftrightarrow & Id_E \circ A = B^{-1} \circ B' \circ A' \\
& \Leftrightarrow&   B^{-1} \circ B' = A \circ A'^{-1} \in PDO^*(M,E) \\
& \Rightarrow &  B^{-1} \circ B' \in DO^{0,*}(M,E) \\
& \Leftrightarrow & \pi(B^{-1} \circ B') = Id_M \\
& \Leftrightarrow & \pi(B)= \pi( B') \end{eqnarray*}
\qed

The next lemma is  obvious:
\begin{Lemma}
$FIO_{Diff}^*(M,E)$ is a group.
\end{Lemma}

\begin{Lemma}
$Ker(\tilde{\pi}) = PDO^*(M,E)$
\end{Lemma} 

\noindent
\textbf{Proof.}
Let $B\circ A \in FIO^*(M,E)$ such that $$\tilde{\pi}(B\circ A)= \pi(B) = Id_M. $$
Then $B \in DO^{0,*}(M,E)$ and $B \circ A \in PDO^*(M,E).$ \qed 

These results show the following theorem:

\begin{Theorem}\label{splitfiodiff}
There is a short exact sequence of groups : 
$$ 0 \rightarrow PDO^*(M,E) \rightarrow FIO^*_{Diff}(M,E) \rightarrow Diff(M) \rightarrow 0$$
and, if $H\subset PDO^*(M,E)$ is a regular Fr\'echet or Fr\"olicher Lie group of operators that contains the gauge group of $E,$
if $K$ is a regular Fr\'echet or Fr\"olicher Lie subgroup of $Diff(M)$ such that there exists a local section  $K \rightarrow Aut(E),$
the subgroup $G = K \circ H$ of $FIO^*_{Diff}(M,E)$ is a regular Fr\'echet Lie group from the short exact sequence:
$$ 0 \rightarrow H \rightarrow G \rightarrow K \rightarrow 0.$$
\end{Theorem}
 \subsection{Diffeomorphisms and kernel operators}
 Let $g \in Diff(M).$ Then a straightforward computation on local coordinates shows that the kernel of $T_g$ is $$K_g = \delta(g(x),y)$$
 where $\delta$ is the Dirac $\delta-$function. These operators also read locally as $$T_g(f) = \int_M e^{ig(x).\xi} \hat{f}(\xi) d \xi$$
on the same system of local trivializations used in \cite{Gil}, p.30-40.

\subsection{Renormalized traces} \label{s3}
 
$E$ is equipped this an Hermitian products $<.,.>$,
which induces the following $L^2$-inner product on sections of $E$: 
$$ \forall u,v \in C^\infty(S^1,E), \quad (u,v)_{L^2} = \int_{S^1} <u(x),v(x)> dx, $$
where $dx$ is the Riemannian volume. 
\begin{Definition} \cite{PayBook,Scott}
$Q$ is a \textbf{weight} of order $s>0$ on $E$ if and only if $Q$ is a classical, elliptic,
 admissible pseudo-differential operator acting on 
smooth sections of $E$, with an admissible spectrum.
\end{Definition}
Recall that, under these assumptions, the weight $Q$ has a real discrete spectrum, and that 
all its eigenspaces are finite dimensional. 
For such a weight $Q$ of order $q$, one can define the complex 
powers of $Q$ \cite{Se}, 
see e.g. \cite{CDMP} for a fast overview of technicalities. 
The powers $Q^{-s}$ of the weight $Q$ 
are defined for $Re(s) > 0$ using with a contour integral,
$$ Q^{-s} = \int_\Gamma \lambda^s(Q- \lambda Id)^{-1} d\lambda,$$
where $\Gamma$ is an ``angular'' contour around the spectrum of $Q.$
Let $A$ be a log-polyhomogeneous pseudo-differential operator.  The map
$\zeta(A,Q,s) = s\in \C \mapsto \hbox{tr} \left( AQ^{-s} \right)\in \C$ , 
defined for $Re(s)$ large, extends
on $\C$ to a meromorphic function with a pole of order $q+1$ at
$0$ (\cite{Le}). When $A$ is classical, 
$\zeta(A,Q,.)$ has a simple pole at $0$  
with residue ${1 \over q} \res A$, where $\res$ is the Wodzicki
residue (\cite{W}, see also \cite{Ka}). Notice that the Wodzicki residue extends the Adler trace \cite{Adl} on formal symbols.
Following \cite{Le}, we define the renormalized trace, see e.g. \cite{CDMP}, \cite{Pay} for the renormalized trace of 
classical operators.

\begin{Definition} \label{d6} $tr^Q A = lim_{z \rightarrow 0} (\hbox{tr} (AQ^{-z})
- {1 \over qz} res A)$.
\end{Definition}

On the other hand, the operator $e^{-tQ}$ is a smoothing operator for each $t>0,$ which shows that $trAe^{-tQ}$ is well-defined and finite for $t>0.$ 
From the function $t\mapsto trAe^{-tQ}, $ we recover the function   $z \mapsto \hbox{tr} (AQ^{-z})$ by the Mellin transform (see e.g. \cite{Scott}, pp. 115-116), which shows the following lemma: 

\begin{Lemma}\label{EQ}
Let $A,A'$ be classical pseudo-differential operators, let $Q,Q'$ be weights.
$$ \forall t>0, trAe^{-tQ}= trA'e^{-tQ'} \Rightarrow \left\{ \begin{array}{ccc} tr^Q(A)&=& tr^{Q'}(A') \\
res(A) & = & res(A')\end{array} \right.$$
\end{Lemma}

If $A$ is trace class, $\hbox{tr}^Q{(A)}=\hbox{tr}{(A)}$.
The functional $\hbox{tr}^Q$ is of course not a trace on $Cl(M,E)$.
Notice also that, if $A$ and $Q$ are pseudo-differential operators
acting on sections on a real vector bundle $E$, they also act on
$E \otimes \C$. The Wodzicki residue res and the renormalized traces
$\hbox{tr}^Q$ have to be understood as functional defined on
pseudo-differential operators acting on $E \otimes \C$.
In order to compute $\tr^Q[A,B]$ and to differentiate $\tr^Q A$,
in the topology of classical pseudo-differential operators, we
need the following (\cite{CDMP}, see also \cite{MN} for the first
point):

\begin{Proposition} \label{p6}

(i)  Given two (classical) pseudo-differential operators A and B,
given a weight Q,
\begin{equation}\label{crochet}  \tr^Q[A,B] = -{1 \over q} \res (A[B,\log Q]). \end{equation}
(ii) Given  a differentiable family $A_t$ of pseudo-differential
operators, given a differentiable family $Q_t$ of weights of
constant order q,

\begin{equation}\label{deriv} {d \over dt} \left(tr^{Q_t}A_t\right) = tr^{Q_t} \left({d \over dt}A_t\right) -
{1 \over q} \res \left( A_t ({d \over dt}\log Q_t) \right).
\end{equation}

\end{Proposition}

The following "covariance" property of $\hbox{tr}^Q$ (\cite{CDMP}, \cite{Pay})
will be useful to define renormalized traces on bundles of operators,

\begin{Proposition} \label{p7}

Under the previous notations, if C is a classical elliptic
injective operator of order 0, $tr^{C^{-1}QC}\left( C^{-1}AC
\right)$ is well-defined and equals $\tr^QA$.

\end{Proposition}

We moreover have specific properties for weighted traces of a more
restricted class of pseudo-differential operators (see
\cite{KV1},\cite{KV2},\cite{CDMP}), called odd class
pseudo-differential operators following \cite{KV1},\cite{KV2} :

\begin{Definition} \label{d7}

A classical pseudo-differential operator $A$ is called odd class
if and only if

$$ \forall n \in \Z, \forall (x,\xi) \in T^*M, \sigma_n(A) 
(x,-\xi) = (-1)^n  \sigma_n(A) (x,\xi).$$

We note this class $Cl_{odd}.$
\end{Definition}

Such a definition is consistent for pseudo-differential operators
on smooth sections of vector bundles, and applying the local
formula for Wodzicki residue, one can prove \cite{CDMP}:

\begin{Proposition} \label{p8}

If $M$ is an odd dimensional manifold, $A$ and $Q$ lie in the odd
class, then $f(s)=tr(AQ^{-s})$ has no pole at $s=0$. Moreover, if
A and B are odd class pseudo-differential operators, $\tr^Q \left(
[A,B] \right) =0$ and $\tr^QA$ does not depend on $Q.$

\end{Proposition}

This trace was first defined in the papers \cite{KV1} and \cite{KV2} by Kontesevich and Vishik. We remark that it is in particular a trace on $DO(M,E)$
when $M$ is odd-dimensional.

Let us now describe a class of operators which is, in some sense, complementary to odd class:
\begin{Definition} 

A classical pseudo-differential operator $A$ is called even class
if and only if

$$ \forall n \in \Z, \forall (x,\xi) \in T^*M, \sigma_n(A) 
(x,-\xi) = (-1)^{n+1}  \sigma_n(A) (x,\xi).$$
 We note this class $Cl_{even}.$
\end{Definition}

Very easy properties are the following:

\begin{Proposition}

$ Cl_{even} \circ Cl_{odd} = Cl_{odd} \circ Cl_{even} = Cl_{even}$
and 

$Cl_{even} \circ Cl_{even} = Cl_{odd} \circ Cl_{odd} = Cl_{odd} .$
\end{Proposition}

Now, following \cite{Ma2006}, we explore properties of $\tr^Q$ on Lie brackets.

\begin{Definition} \label{d8}

Let E be a vector bundle over M, Q a weight and $a \in \Z$. We define :
$$ \mathcal{A}^Q_a=\{B \in Cl(M,E); [B,\log Q] \in Cl^a(M,E)\}.$$

\end{Definition}

\begin{Theorem} \label{t1} \cite{Ma2006}

\begin{item}

(i) $\mathcal{A}^Q_a \cap Cl^0(M,E) $ is an subalgebra of $Cl(M,E)$
with unit.
\end{item}

\begin{item}

(ii) Let $B \in Ell^*(M,E)$, $B^{-1}\mathcal{A}^Q_aB = A^{B^{-1}QB}_a.$

\end{item}

\begin{item}

(iii) Let $A\in Cl^b(M,E)$, and $B \in \mathcal{A}^Q_{-dimM-b-1}$,
then $\tr^Q[A,B]=0.$
\end{item}

\begin{item}

(iv)  For $a< -{dimM \over 2}$, $\mathcal{A}^Q_a \cap Cl^{-dimM \over 2}(M,E) $
is an algebra on which the renormalized trace is a trace (i.e. vanishes on the brackets).

\end{item}

\end{Theorem}

We now produce non trivial examples of operators that
are in $\mathcal{A}^Q_{a}$ when Q is scalar, and secondly we give a
formula for some non vanishing renormalized traces of a bracket.

\begin{Lemma} \label{l1} Let Q be a weight on $C^\infty_0(M, V)$ and
let B be a classical pseudo-differential operator of order $b$. If
$B$ or $Q$ is scalar, then $[B, \log Q]$ is a classical
pseudo-differential operator of order $b-1$.

\end{Lemma}

\begin{Proposition} \label{p9}

Let Q be a scalar weight on $C^\infty_0(M, V)$. Then
$$Cl^{a+1}(M,V) \subset \mathcal{A}^Q_{a}.$$
Consequently,

\begin{item}

(i) if $ord(A)+ord(B)=-dimM,$ $\tr^Q [A,B]=0$.

\end{item}

\begin{item}

(ii) when $M=S^1$, if A and B are classical pseudo-differential
operators, if A is compact and B is of order 0, $\tr^Q[A,B]=0$.

\end{item}
\end{Proposition}

\begin{Lemma} \label{l3}

Let Q be a scalar weight on  $ C^\infty_0(M,V) $, and A, B two
pseudo-differential operators of orders a and b on $
C^\infty_0(M,V) $, such that $a+b=-m+1$ (m = dim M). Then
$$ \tr^Q[A,B]=-{1 \over q} \res\left( A[B,\log Q] \right)=
-{1 \over q(2\pi )^n} \int _M \int _{|\xi |=1} tr (\sigma _a( A
)\sigma_{b-1}([B, \log Q])). $$

\end{Lemma}

Let us now explore the action of $Diff(M)$ and of $Aut(E)$ on $tr^Q(A).$For this, we get: 

\begin{Lemma}
Let $a \in \Z.$ Let $A \in Cl^a(M,E)$ and let $Q$ be a weight on $E.$ Let $B$ be an operator on $C^{\infty}(M,E)$ such that 
\begin{enumerate}
\item \label{a} $Ad_B(Cl^a(M,E)) \subset Cl^a(M,E)$
\item \label{b} $Ad_BQ$ is a weight of the same order as $Q$
\end{enumerate} 
Then 
\begin{itemize}
\item $res(Ad_BA) = res(A)$
\item $ tr^{Ad_BQ}(Ad_BA) = tr^Q(A).$
\end{itemize}
The properties \ref{a},\ref{b} are true in particular for operators $B\in Aut(E).$
\end{Lemma}

\noindent
\textbf{Proof.}

Let $Q$ be a weight on $C^\infty(M,E)$ and let $A \in Cl(M,E).$ Let $B \in Aut(E).$ Let $s\in \R_+^{*}$ then $Ae^{-sQ}$ is trace class. By \cite{Gil}, we know that $Ad_BA$ (resp. $Ad_BQ$) is a classical pseudo-differential operator of the same order (resp. a weight of the same order). 
Then, since $e^{-\frac{s}{2}Q}$ is smoothing, $Ad_{B}(Ae^{-sQ},$  $BAe^{-\frac{s}{2}Q}$ and $e^{-\frac{s}{2}Q}B^{-1}$ are smoothing, and the following computations are fully justified:

\begin{eqnarray*}
tr\left(Ad_{B}(Ae^{-sQ})\right)& = &  tr\left(\left( BAe^{-\frac{s}{2}Q} \right)\left( e^{-\frac{s}{2}Q}B^{-1} \right)\right) \\
& = & tr\left(\left( e^{-\frac{s}{2}Q}B^{-1} \right)\left( BAe^{-\frac{s}{2}Q} \right)\right) \\
& = & tr\left( e^{-\frac{s}{2}Q}Ae^{-\frac{s}{2}Q} \right) \\
& = & tr\left( Ae^{-sQ} \right) \\
\end{eqnarray*}
So that, by Lemma \ref{EQ},
we get the announced property. \qed

\section{Splittings on the set of $S^1-$Fourier integral operators}

\subsection{The group $O(2)$ and the diffeomorphism group $Diff(S^1)$}
Let us consider the $SO(2)=U(1)$-action on $S^1=\R/ \mathbb{Z}$ given by $(e^{2i\pi\theta},x)\mapsto x + \theta.$ This group acts on $C^\infty$ by 
$ (e^{2i\pi\theta}, f) \mapsto f(x+\theta)$ and we have
\begin{eqnarray*} f (x + \theta) & = & \int e^{-i(x+\theta).\xi} \hat{f}(\xi) d\xi \\
				& = & \int e^{-i(x.\xi + \theta.\xi} \hat{f}(\xi) d\xi
\end{eqnarray*}
The term $e^{-i\theta.\xi}$ is oscillating in $\xi$ and does not satisfies 
the estimates on the derivatives of symbols. So that, this operator 
is not a pseudo-differential operator but has obviously the form of a 
Fourier integral operator. 
The same is for the reflection $x \mapsto 1-x$ which corresponds to the conjugate
transformation $z \mapsto \bar{z}$ when representing $S^1$ as the set
of complex numbers $z$ such that $|z|=1.$
This is a spacial case of the properties already stated for a general manifold $M$ given $g \in Diff(S^1),$ $g$ acts on $C^\infty$ by right 
composition of the inverse, namely, for $f \in C^\infty,$
\begin{eqnarray*} g.f (x)& = & f \circ g(x)\\
			& = & \int e^{-ig(x).\xi} 
				\hat{f}(\xi) d\xi, \end{eqnarray*}
which is also obviously a Fourier-integral operator,
 and the kernel of this operator is $$K_g (x,y) = \delta\left(y , g(x)\right)$$
where $\delta$ is the Dirac $\delta$-function. This is the construction already used in the proof of Theorem \ref{FIOdiff}.

\subsection{$\epsilon(D),$ its formal symbol and the splitting of $\F PDO$} 
The operator $D = {-i} D_x$ splits $C^\infty(S^1, \C^k)$ into three spaces :

- its kernel $E_0$, made of constant maps

- $E_+$, the vector space spanned by eigenvectors related to positive eigenvalues

- $E_-$, the vector space spanned by eigenvectors related to negative eigenvalues.

\noindent
The following elementary result will be useful for the sequel, 
see \cite{Ma2003} for the proof, and e.g. \cite{Ma2006,Ma2008}:  
\begin{Lemma} \label{l1'}

(i) $\sigma(D) = {\xi }$

(ii) $\sigma(|D|) = {|\xi| }$ where $|D| =
 \sigma \left(\int_\Gamma \lambda^{1/2} (\Delta - \lambda Id)^{-1}d\lambda\right)$, 
with $\Delta = -D_x^2$. 

(iii)  $\sigma(D|D|^{-1}) = {\xi \over |\xi|}$, where $ D|D|^{-1} = |D|^{-1}D$ is the sign of D, since $|D|_{|E_0}=Id_{E_0}.$

(iv)  Let $p_{E_+}$ (resp. $p_{E_-}$) be the projection on $E_+$ (resp. $E_-$), then 
$\sigma(p_{E_+}) ={1 \over 2}(Id + {\xi \over |\xi|})$ and $\sigma(p_{E_-}) = {1 \over 2}(Id - {\xi \over |\xi|})$.

\end{Lemma}

Let us now define two ideals of the algebra $\mathcal{F}PDO$, 
that we call $\mathcal{F}PDO_+$ and $\mathcal{F}PDO_-$, 
such that $\mathcal{F}PDO = \mathcal{F}PDO_+ \oplus \mathcal{F}PDO_-$. 
This decomposition is implicit in \cite{Ka}, section 4.4., p. 216,
for classical pseudo-differential operators 
and we furnish the explicit description 
given in \cite{Ma2003}, extended to the whole algebra of 
(maybe non formal, non classical) pseudo-differential symbols here.

\begin{Definition}

Let $\sigma$ be a symbol (maybe non formal). Then, we define, for $\xi \in T^*S^1 - S^1$, 
$$ \sigma_+(\xi) = \left\{ 
\begin{array}{ll}
\sigma(\xi) & \hbox{ if $ \xi > 0$} \\
0 & \hbox{ if $ \xi < 0$} \\
\end{array}
\right. \hbox{ and }
 \sigma_-(\xi) = \left\{ 
\begin{array}{ll}
0 & \hbox{ if $ \xi > 0$} \\
\sigma(\xi) & \hbox{ if $ \xi < 0$} . \\
\end{array}
\right.$$
At the level of formal symbols, we also define the projections:  $p_+(\sigma) = \sigma_+$ and $p_-(\sigma) = \sigma_-$ .
\end{Definition}
The maps 
$ p_+ : \mathcal{F}PDO(S^1,\C^k) \rightarrow \mathcal{F}PDO(S^1,\C^k) $ 
{ and } $p_- : \mathcal{F}PDO(S^1,\C^k) \rightarrow \mathcal{F}PDO(S^1,\C^k)$ 
are clearly  algebra morphisms 
that leave the order invariant and are also projections 
(since multiplication on formal symbols is expressed 
in terms of pointwise multiplication of tensors). 

\begin{Definition} We define
$  \mathcal{F}PDO_+(S^1,\C^k) = Im(p_+) = Ker(p_-)$
and $  \mathcal{F}PDO_-(S^1,\C^k) = Im(p_-) = Ker(p_+).$ \end{Definition}
Since $p_+$ is a projection,  we have the splitting
$$ \mathcal{F}PDO(S^1,\C^k) = \mathcal{F}PDO_+(S^1,\C^k) \oplus \mathcal{F}PDO_-(S^1,\C^k) .$$
Let us give another characterization of $p_+$ and $p_-$. 
Looking  more precisely at the formal symbols of $p_{E_+}$ and $p_{E_-}$ 
computed in Lemma \ref{l1'}, we observe that 
$$ \sigma( p_{E_+}) = \left\{ \begin{array}{ll}
1 & \hbox{if }\xi > 0 \\
0 & \hbox{if }\xi < 0 \\
\end{array} \right. \hbox{ and }
 \sigma( p_{E_-}) = \left\{ \begin{array}{ll}
0 & \hbox{if }\xi > 0 \\
1 & \hbox{if }\xi < 0 \\
\end{array} \right. . $$
In particular, we have that $D^\alpha_x\sigma( p_{E_+}),$ 
$D^\alpha_\xi\sigma( p_{E_+}),$ $D^\alpha_x\sigma( p_{E_-}),$ $D^\alpha_\xi\sigma( p_{E_-})$ 
vanish for $\alpha > 0$. 
From this, we have the following result:

\begin{Proposition} \label{pag} \cite{Ma2003}
Let $a \in \mathcal{F}PDO(S^1,\C^k).$ 
$ p_+(a) =  \sigma( p_{E_+}) \circ a = a \circ \sigma( p_{E_+})$ and 
$  p_-(a) =  \sigma( p_{E_-}) \circ a = a \circ \sigma( p_{E_-}).$
\end{Proposition}

\subsection{The case of non trivial (real) vector bundle over $S^1$}

Let $\pi : E \rightarrow S^1$ be a non trivial real vector bundle over $S^1$ of rank $k.$ Its bundle of frames is a $Gl(\R^k)-$ principal bundle, which means the following (see e.g. \cite{KN}):
\begin{Lemma}
Let $\varphi_1 : ]a;b[ \times \R^k \rightarrow E$ and $\varphi_2 : ]a';b'[ \times \R^k \rightarrow E$ be two local trivializations of $E$. Let $\mathcal{D} = \pi( \varphi_1 ( ]a;b[ \times \R^k) \cap \varphi_2 ( ]a';b'[ \times \R^k)),$ let $\mathcal{D}_1 = \varphi_1^{-1}(\mathcal{D}),$ and let $\mathcal{D}_2 = \varphi_2^{-1}(\mathcal{D}). $ Then $$ \varphi_2^{-1} \circ \varphi_1 : \mathcal{D}_1 \times \R^k \rightarrow \mathcal{D}_2 \times \R^k$$
reads as $$ \varphi_2^{-1} \circ \varphi_1  = \gamma \times M$$
where $ \gamma$ is a smooth diffeomorphism from $\mathcal{D}_1$ to $\mathcal{D}_2,$ and where $M \in C^\infty(\mathcal{D}_1, Gl(\R^k)).$
\end{Lemma} 

Let us now turn to symbols of pseudo-differential operators acting on smooth sections of $E.$ We first assume that we work with a system of local trivializations such that the diffeomorphisms $\gamma$ are translations, and let us now look at the transformations of the symbols read on local trivializations. Under these assumptions, and with the notations of the previous lemma, a formal symbol $\sigma_1$ read on $D_1$ reads on $D_2$ as 
$$ \sigma_2(\gamma(x),\xi) = M(x) \sigma_1(x,\xi) M(x)^{-1}.$$

\begin{Proposition}
Let $\nabla$ be a Riemannian covariant derivative on the bundle $E \rightarrow S^1$ and let $\nabla \over dt$ be the associated first order differential operator, given by the covariant derivative evaluated at the unit vector field over $S^1.$  We modify the operator $\nabla \over dt$ into an injective operator $D = {\nabla \over dt} + p_{ker {\nabla \over dt}}$, where $p_{ker {\nabla \over dt}}$ is the $L^2$ orthogonal projection on $ker {\nabla \over dt} \subset C^\infty(S^1,E) \subset L^2(S^1,E),$ and we set $$\epsilon(\nabla) =  D \circ \left| D \right|^{-1}.$$
Then the formal symbol of  $\epsilon(\nabla)$ is $i\xi \over |\xi |.$
\end{Proposition}

\noindent
\textbf{Proof.}
Let us use the holonomy trivialization over an interval $I.$ In this trivialization,  $${\nabla \over dt} = {d \over dt} $$
and hence the formal symbol of ${\nabla \over dt}$ reads as $i\xi.$ Calculating exclusively on the algebra of formal operators on which composition and inversion governed by local formulas, we get $\sigma(|D|) = |\xi|$ and, by the same arguments as those of \cite{Ma2003}, we get the result. \qed

\begin{Proposition}
For each $A \in PDO(S^1,E),$  $[A,\epsilon(\nabla)] \in PDO^{-\infty}(S^1;E).$
\end{Proposition}

\textbf{Proof.}
We remark that, for any multiindex $\alpha$ such that $|\alpha| >0,$ $D_x^\alpha \sigma(\epsilon(\nabla)) =0$ and $D^\alpha_\xi \sigma(\epsilon(\nabla))=0.$ Hence, in $\F PDO (S^1,E),$ $$\sigma ([A,\epsilon(\nabla)]) = [\sigma (
A),\sigma(\epsilon(\nabla))]=0$$
so that $ [A,\epsilon(\nabla)] \in PDO^{-\infty}(S^1,E).$ \qed

\subsection{The splitting read on the phase function} 
The fiber bundle $T^*S^1 - S^1$ has two connected components and the 
phase function is positively homogeneous, so that
we can make the same procedure as in the case of the symbols. 
But we remark that we can split
$$ \phi = \phi_+ + \phi_-$$
where  $\phi_+ = 0$ if $\xi <0$ and $\phi_- = 0$ if $\xi >0.$ Unfortunately, $\phi_+$ and $\phi_-$ 
are \textbf{not} phase functions of Fourier integral operators because there are some points where ${\partial^2 \phi_+ \over \partial_x \partial_\xi } = 0$ or ${\partial^2 \phi_- \over \partial_x \partial_\xi } = 0.$   
However, we can have the following identities:
\begin{eqnarray*} \int_\R e^{i\phi(x,\xi)} \sigma(x,\xi) \hat{f}(\xi) d\xi & = & 
\int_{\xi >0} e^{i\phi(x,\xi)} \sigma(x,\xi) \hat{f}(\xi) d\xi + \int_{\xi < 0} e^{i\phi(x,\xi)} \sigma(x,\xi) \hat{f}(\xi) d\xi \\
& = & 
\int_{\xi >0} e^{i\phi_+(x,\xi)} \sigma(x,\xi) \hat{f}(\xi) d\xi + \int_{\xi < 0} e^{i\phi_-(x,\xi)} \sigma(x,\xi) \hat{f}(\xi) d\xi \\
& = & 
\int_{\R} e^{i\phi_+(x,\xi)} \sigma_+(x,\xi) \hat{f}(\xi) d\xi + \int_{\R} e^{i\phi_-(x,\xi)} \sigma_-(x,\xi) \hat{f}(\xi) d\xi  \\
& = & 
\int_{\R} e^{i\phi(x,\xi)} \sigma_+(x,\xi) \hat{f}(\xi) d\xi + \int_{\R} e^{i(x,\xi)} \sigma_-(x,\xi) \hat{f}(\xi) d\xi\end{eqnarray*}

\subsection{The Schwinger cocycle on $PDO(S^1,E)$ when $E$ is a real vector bundle.} \label{sect2}

The main result of \cite{Ma2006,Ma2008} are now analyzed from the viewpoint of operators acting on smooth sections of real vector bundles. Here, $\epsilon(\nabla)$ is not a sign operator, but an operator such that $\epsilon(\nabla)^2 = -Id$ up to a smoothing operator. :  

\begin{Theorem} \label{th1}
For any $A \in PDO(S^1,E)$, $[A,\epsilon(\nabla)] \in PDO^{-\infty}(S^1,E).$ Consequently, 
$$ c_s^\nabla : A,B \in PDO(S^1,E) \mapsto {1 \over 2}\tr \left( \epsilon(\nabla)[\epsilon(\nabla),A][\epsilon(\nabla),B]  \right) $$
is a well-defined $\mathbb{R}$-valued 2-cocycle on $PDO(S^1, E).$ Moreover, $c_s^\nabla$ is non trivial on any Lie algebra $\mathcal A$ such that $C^\infty(S^1,\R)\subset \mathcal{A} \subset PDO(S^1,E).$
\end{Theorem}
Notice that $C^\infty(S^1,\R)$ is understood as an algebra acting on $C^\infty(S^1, E)$ by scalar multiplication fiberwise. The proof follows the same arguments as in \cite{Ma2008}.

\vskip 12pt
\noindent
\textbf{Proof.} First, $c_s^\nabla$ is the trace of operators acting on a real Hilbert space. so that, it is real valued.
Since \cite{PS}, see e.g. \cite{Ma2006}, if $c_s^\nabla$ was trivial on Hoschild cohomology, there would have a 1-form $\nu : \mathcal{A} \rightarrow \R$ such that $$c_s^D = \nu([.,.]),$$ and hence it would be true on $C^\infty(S^1,\R)$ which is a commutative algebra. Hence, since $c_s^\nabla \neq 0$ on $C^\infty(S^1,\R),$ it is non trivial on $A.$ \qed
\section{Sets of Fourier Integral operators}
\subsection{The set $FIO(S^1,E)$}

Here, for the definitions, $\eps = \epsilon(D)$ or $\eps=\eps(\nabla),$ depending on the fact that $E$ is a complex or a real vector bundle.
Let us now define
$$ FIO_{res}(S^1;E) = \{ A \in FIO(S^1,E) \hbox{ such that } [A;\epsilon]\in PDO^{-\infty}(S^1,E)\}.$$

\begin{Proposition}
$FIO_{res}(S^1,E)$ is a set, stable under composition, with unit element.
\end{Proposition}

\vskip 12pt
\noindent
\textbf{Proof.} $FIO(S^1,E)$ is stable under composition \cite{Horm}. 
Since $Cl^0(S^1,E)$ is contained in $FIO_{res}(S^1,E)$ by Theorem \ref{th1} so that 
$FIO_{res}(S^1,E)$ contains the identity map. 

Let $A,B \in FIO_{res}(S^1,E),$ 
\begin{eqnarray*}
[AB,\epsilon] & = & A[B;\epsilon] + [A;\epsilon]B
\end{eqnarray*}
Since $[A,\epsilon]$ and $[B;\epsilon]$ are smoothing, we get that $[AB, \epsilon]$ is smoothing. \quad \qed
 
We use the natural notations, $$FIO_{res}^0=FIO^0 \cap FIO_{res}.$$

We shall note by $ FIO_{res}^*(S^1,E)$ the group of units of this set, 
and by  $ FIO_{res}^{0,*}(S^1,E)$ the group of units of the set
$ FIO_{res}^0(S^1,E). $

\begin{Proposition}
$ FIO_{res}^*(S^1,E) = FIO^{*}(S^1,E) \cap FIO_{res}(S^1,E)$ 
and $ FIO_{res}^{0,*}(S^1,E) = FIO^{0,*}(S^1,E) \cap FIO_{res}(S^1,E)$
\end{Proposition}

\noindent
\textbf{Proof.}
We already have trivially $ FIO_{res}^*(S^1,E) \subset (S^1,E)\cap FIO_{res}(S^1,E). $
Let $A \in FIO^{*}(S^1,E) \cap FIO_{res}(S^1,E).$ 
We have to check that $A^{-1} \in  FIO_{res}(S^1,E).$
\begin{eqnarray*}
A[A^{-1},\eps] & = & [A A^{-1}, \eps] - [A, \eps ]A^{-1} \\
& = & [Id, \eps] - [A, \eps ]A^{-1} \\
& = & - [A, \eps ]A^{-1} \\
& \in & PDO^{-\infty}(S^1,E) \end{eqnarray*}
So that
\begin{eqnarray*}
[A^{-1},\eps] & = & A^{-1} A [A^{-1},\eps] \\
& \in & PDO^{-\infty}(S^1,E) \end{eqnarray*}
The proof is the same for $0-$order operators. \qed

By the way, since $FIO^{0,*}(S^1,E)$ is a "`generalized Lie group"' in the sense of Omori, it is a Fr\"olicher Lie group. By the trace property of Fr\"olicher spaces, using the last proposition,  $FIO_{res}^{0,*}(S^1,E)$ is a Fr\"olicher Lie group \cite{Ma2013}. 

Now, since we have that $$ FIO_{res}^{0,*} \subset GL_{res},$$
the determinant bundle defined over $GL_{res}$ can be pulled-back on $FIO_{res}^{0,*}$. 
The same way, it is shown in \cite{Ma2006-2,Ma2008} that the Schwinger cocycle extends to the Lie algebra 
$PDO^0(S^1,E) + PDO^1(S^1,\C)\otimes Id_E.$

\subsection{Yet some subgroups of $ FIO_{res}^*(S^1,E)$}

 Let us first gather and reformulate many known results:
\begin{Lemma}
$Diff^+(S^1) \times C^{\infty}(S^1,\C^*) \subset FIO^{0,*}_{res}(S^1,\C).$
\end{Lemma}

\noindent
\textbf{Proof.}
First, we have that $$ C^{\infty}(S^1,\C^*) \subset Cl^{0,*}(S^1,\C)$$
so that 
$$ C^{\infty}(S^1,\C^*) \subset FIO^{0,*}(S^1,\C).$$
Let $g \in Diff^+(S^1).$
Following \cite{PS}, the map 
$f \mapsto |g'|^{1/2}. (f \circ g)$ describes an operator in $U_{res} \subset GL_{res}.$
Since the map $f \mapsto |g'|^{1/2}. f$ is a multiplication operator in $C^{\infty}(S^1,\C^*),$ 
we get that $$f \mapsto  f \circ g = \int e^{-i g(.).\xi} \hat{f}(\xi) d\xi \in GL_{res} \cap FIO^{0,*}(S^1,\C).$$ \qed

\begin{Theorem} 
Assume that $E$ be a trivial vector bundle over $S^1.$ 
Let $\tilde{\pi}$ be the projection $ FIO_{Diff}^{*}(S^1,E) \rightarrow Diff(S^1).$ Then $$\pi^{-1}(Diff_+(S^1)) \subset FIO_{res}(S^1,E).$$
\end{Theorem}
This is a simple consequence of the previous results.

\begin{Theorem}
Assume that $E$ is non trivial and let $\epsilon$ defined as before. Let $\tilde{\pi}$ be the projection $ FIO_{Diff}^{*}(S^1,E) \rightarrow Diff(S^1).$ Then $$FIO^*_{Diff_+}(S^1,E) = \pi^{-1}(Diff_+(S^1)) \subset FIO_{res}(S^1,E),$$
and there is a global smooth section (in the sense of Fr\"olicher spaces, not necessarily in the sense of groups) $$ Diff_+(S^1) \rightarrow FIO_{res}(S^1,E)$$
of the short exact sequence:
$$0 \rightarrow PDO^*(S^1,E) \rightarrow FIO_{Diff}^{*}(S^1,E)\cap FIO_{res}(S^1,E) \rightarrow Diff_+(S^1) \rightarrow 0.$$
\end{Theorem}

\noindent
\textbf{Proof.}
Let  $g\in Diff_+(S^1).$ We fix on $E$ a connection $\nabla$ and we set $n = rank(E).$
Since $Diff_+(S^1)$is the connected component of $Id_{S^1}$ in $Diff(S^1)$, given $\eta$ the unit vector field defined by orientation on $S^1$, we can choose 
a path $$\gamma \in C^\infty([0,1], Diff_+(S^1))\subset C^\infty([0,1] \times S^1, S^1)$$
such that $$\gamma(0) = Id_{S^1}, \gamma(1) = g$$
and  $$\forall x \in S^1, \forall t \in [0;1], (\frac{d \gamma}{dt}(t)(x),\eta(x))_{T_xS^1} >0.$$
This path is unique up to parametrization since we impose also the condition of minimal length.
Let $$H_x = Hol(\gamma(.)(x)) \in Gl(E_x, E_{g(x)})$$
be the induced parallel transport map.
We get, for each $g \in Diff_+(S^1),$
a map $H_g$
which is smooth by the properties of parallel transport, linear on the fibers, invertible, and which projects on $S^1$ to $g.$
Thus, $H_g \in Aut(E),$ and it easy to see that it is a bijection on the collection of smooth trivializations of $E.$ Now,turning to the map $$g \mapsto H_g,$$
is appears as a smooth map $Diff(S^1) \rightarrow Aut(E),$ but it is not a group morphism in any case since $E$ can be non trivial.     We have moreover that 
$$\forall g \in Diff^+(S^1), [\nabla, H_g] = 0$$
since $dim(S^1)=1$ and $H_g$ is a parallel transport map. So that, since $\epsilon$ is derived from ${\nabla \over dt} = \nabla_\eta,$ we get that
$H_g \in GL_{res}.$ Now, an operator in $FIO^*_{Diff_+}(S^1,E)$ reads as $$H_g \circ A,$$ where $A \in PDO^*(S^1,E) \subset GL_{res}.$ Then $H_g \circ A \in FIO_{res}.$\qed

\begin{Theorem}
The group $$FCl^{0,*}_{Diff_+}(S^1,E) = FIO^*_{Diff_+}(S^1,E) \cap FCL^0(S^1,E)$$
is a regular Fr\"olicher Lie group.
\end{Theorem}

\noindent
\textbf{Proof.}
We get the obvious exact sequence of Lie groups:
$$0 \rightarrow Cl^{*,0}(S^1,E) \rightarrow FCl_{Diff}^{*,0}(S^1,E)
\rightarrow Diff_+(S^1) \rightarrow 0.$$
Both $Cl^{*,0}(S^1,E)$ and $Diff_+(S^1)$ are regular, 
and $Aut(E) \subset Cl^{*,0}(S^1,E),$
so that the smooth section $Diff_+(S^1) \rightarrow Aut(E)$ described in the proof of the previous theorem gives the result by Theorem \ref{splitfiodiff}. \qed

Let us now describe a subgroup of $FIO^*_{Diff_+}(S^1,E). $ 
\begin{Definition}
Let $FIO^*_{b,Diff_+}(S^1,E)$ be the space of operators $A \in  FIO^*_{Diff_+}(S^1,E)$ such that 
\begin{enumerate}
\item $\pi(A)$ is a diffeomorphism of $S^1 = \R /\Z$ such that $\pi(A)(0) = 0;$
\item if $u $ is a smooth section of $E$ such that $u(0) = 0,$ then $(Au)(0)=0.$
\end{enumerate}
These operators are called based operators, and the set of sections $u$  of $E$ such that $u(0)=0$ is the space of based sections, noted $C^\infty_b(S^1;E).$ We note by $Diff_{b,+}(S^1)$ the infinite dimensional Lie group of diffeomorphisms $g$ such that $g(0)=0.$ 

\end{Definition}
 
\noindent We recall that $Diff_{b,+}(S^1)$ is a regular Lie subgroup of ${Diff_+}(S^1).$ Its Lie algebra is noted $\mathfrak{diff}_b(S^1).$   
Given $c(t)$ a smooth curve in $FIO^*_{b,Diff_+}(S^1,E),$ starting at $Id_E,$  $\frac{dc(t)}{dt}|_{t=0} = U + X,$
where $X \in \mathfrak{diff}_b(S^1)$ and $ U \in PDO(S^1,E)$ which stabilize $C^\infty_b(S^1;E).$

\begin{Theorem} 
\begin{itemize}
\item Let $G \subset PDO^*(S^1,E)$ be a regular Lie group of based operators, that contains the space of based invertible multiplication operators, with Lie algebra $\mathfrak{g}.$
\item Let $D \subset  Diff_{b,+}(S^1)$ be a regular Lie subgroup of based diffeomorphisms, with regular Lie algebra $\mathfrak{d}.$ 
\end{itemize}
There is a regular Lie group $FG_D \subset FIO^*_{b,Diff_+}(S^1,E)$ for which the following sequence is exact:
$$ 0 \rightarrow G \rightarrow FG_D \rightarrow D \rightarrow O.$$
\end{Theorem}

\noindent
\textbf{Proof.} We consider first the regular Lie group of automorphisms $\pi^{-1}(D) \subset Aut(E).$
Then, with the same arguments, $G$ and generate a group that we note $FG_D,$ and adapting the computations of Lemma \ref{phase}, we obtain the above exact sequence. 
Finally, by Theorem \ref{slice}, $FG_D$ is a regular Lie group. \qed 
\section{Manifolds of embeddings}
\noindent \textbf{Notation :} Let $E \rightarrow M$
be a smooth vector bundle over $M$ with typical fiber $x$.
For $k \in \N^*$,
we denote by

- $E^{\times k}$ the product bundle, of basis $M$, with typical
fiber $F^{\times k}$;

- $\Omega^k(E)$ the space of $k-forms$ on $M$ with values in $E$,
that is, the set of smooth maps $(TM)^{\times k} \rightarrow E$ that are fiberwise $k$-linear
and skew-symmetric $(T_xM)^{\times k} \rightarrow E_x$ for any $x \in M$. If $E = M \times F$,
we note $\Omega^k(M,F)$ the space of $k$-forms instead of $\Omega^k(E)$.

\vskip 5pt

Let $M$ be a compact manifold without boundary; let $N$ be a
Riemannian manifold, equipped with the metric $(.,.)$. Let
$Emb(M,N)$ be the manifold of smooth embeddings $M\rightarrow N$.

\subsection{$Emb(M,N)$ as a principal bundle} 
The group of diffeomorphisms of $M$, $Diff(M)$,
acts smoothly and on the right on $Emb(M,N)$, by composition. Moreover,
$$B(M,N)=Emb(M,N)/ Diff(M)$$ is a smooth manifold \cite{KM}, and $\pi: Emb(M,N)\rightarrow B(M,N)$
is a principal bundle with structure group $Diff(M)$ (see
\cite{KM}). Then, $g \in Emb(M,N)$ is in the $Diff(M)-$orbit of $f$ if and only if
$g(M)=f(M)$. 
Let us now precise the vertical tangent space and a normal vector
space of the orbits of $Diff(M)$ on $Emb(M,N)$. $T_fPEmb(M,N)$,
the tangent space at $f$, is identified with the space of smooth
sections of $f^*TN$, which is the pull-back of $TN$ by $f$.
$VT_fP$, the vertical tangent space at $f$ is the space of smooth
sections of $Tf(M)$. Let $\normal_f$ be the normal space to
$f(M)$ with respect to the metric $(.,.)$ on $N$. For any $x \in M$,
$T_{f(x)}N = T_{f(x)}f(M) \oplus \normal f(M)$. Hence, denoting $f* \normal_f$ the
pull back of $\normal _f $ by $f$, we have that
$$ C^\infty(f^*TN) = C^\infty (TM) \oplus f^*\normal_f.$$
Moreover, for any volume form $dx$ on $M$,
if $$ <.,.> : X,Y \in C^\infty(f^*TN) \mapsto <X,Y> = \int_M (X(x),Y(x)) dx$$ is a
$L^2$-inner product on $C^\infty(f^*TN)$, this splitting is orthogonal for $<.,.>$.
We get here a fundamental difference between the inclusion $Emb(M,N) \subset C^\infty(M,N),$
where the model space of the type $C^\infty(f^*TN),$ and $Emb(M,N)$ as a $Diff(M)-$ principal bundle: 
sections of the vertical
tangent vector bundle read as order 1 differential operators, 
where as the operators acting on the normal vector bundle reads
as $0-$order differential operators, just like the structure group of $TC^\infty(M,N).$ 
To be more precise, let $X \in C^\infty(f^*TN)$ and let $p: f^*TN \rightarrow Tf(M)$ be the orthogonal projection. 
The vector field $p(X) \in C^\infty(Tf(M))$ is seen as a differential operator acting on smooth functions 
$f(M)\sim M \rightarrow \R,$ and the normal component $(Id-p)(X)$ is a smooth section on $\normal_f.$  In the sequel we shall note $$\normal = \coprod_{f \in Emb(M,N)} \normal_f.$$

We turn now to local trivializations. 
Let $f \in C^\infty_b(M,N)$. We define the map $Exp_f :
C^\infty_0(M,f^*TN)\rightarrow C^\infty_b(M,N)$ defined by
$Exp_f(v)=exp_{f(.)}v(.)$ where $exp$ is the exponential map on
$N.$ Then $Exp_f$ is a smooth local diffeomorphism. Restricting
$Exp_f$ to a $C^\infty$ - neighborhood $\tilde U_f$ of the 0-section
of $f^*TN$, we define a diffeomorphism, setting

$$ (Exp_f)_{|\tilde U_f} : \tilde U_f
\rightarrow V_f=Exp_f( \tilde U_f) \subset C^\infty_b(M,N) .$$

Then, setting $U_f = I_f^{-1} \tilde U_f$, we can define a chart
$\Xi^f$ on $V_f$ by:

$$ \Xi^f(g) = (I_f^{-1} \circ (Exp_f)_{|\tilde U_f}^{-1})(g) \in U_f \subset
C^\infty_b(M,E) .$$

Given $f,g$ in $C^\infty_b(M,N)$ such that $V_{f,g}= V_f \cap V_g
\neq \not0,$ we compute the changes of charts $\Xi^{f,g}$ from
$U^f_{f,g} = \Xi^fV_{f,g}$ to $U^g_{f,g} = \Xi^gV_{f,g}$. Let $u
\in U^f_{f,g}$, $v=(\Xi^f)^{-1}(u) \in V_{f,g}$.

$$\Xi^{f,g}(u)= \Xi^g \circ (\Xi^f)^{-1} (u)=(I_g^{-1} \circ (Exp_g)^{-1} \circ
Exp_f \circ I_f)(u).$$

 Since, $\forall x \in M$, the transition maps

$$\Xi^{f,g}(u)(x)=(I_g^{-1} \circ (exp_{g(x)})^{-1} \circ exp_{f(x)} \circ
I_f)(u(x))$$ are smooth,  $(V_f, \Xi^f, U_f)_{f \in
C^\infty_b(M,N)}$ is a smooth atlas on $C^\infty_b(M,N)$.
Moreover, let $w \in C^\infty_0(M,E)$, setting
$v=(\Xi^f)^{-1}(u)$, the evaluation of the differential at $x\in
M$ reads :
$$ D_u\Xi^{f,g}(w)(x) =
(I_g^{-1} \circ D_{v(x)}(exp_{g(x)})^{-1} \circ
D_{u(x)}(exp_{f(x)} \circ I_f))(w(x)).$$

Hence, for $u \in C^\infty$,  $D_u\Xi^{f,g}$ is a multiplication
operator acting on smooth sections of E for any isomorphism $I_f$
and $I_g$ we can choose. Since $I_f$ and $I_g$ are fixed, the
family $ u \mapsto  D_u\Xi^{f,g} $ is a smooth family of 0- order
differential operators; this construction is described carefully in \cite{Ee}.
Now, let $f \in Emb(M,N)$ and let us consider the map $$ \Phi^{U,f}: (f,v, X) \in TU \sim (1-p)TU \oplus pTU \mapsto \Xi^f(v). exp_{Diff(M)}(X) \in Emb(M,N) $$ where $p$ is the orthogonat projection on $f_*TS^1.$
This map gives a local (fiberwise) trivialization of the principal bundles $Emb(M,N) \rightarrow B(M,N)$ following \cite{HV,KM,Mo}, and we see that the changes of local trivializations have $Aut(\mathcal{N})$ as a structure group.

If $M$ is oriented, we note by $Diff_+(M)$ the group of orientation preserving diffeomorphisms and we have the following trivial lemma: 

\begin{Lemma}
$$\frac{Diff(M)}{Diff_+(M)} = \mathbb{Z}_2.$$
\end{Lemma}
Then, defining $$B_+(M,N) = \frac{Emb(M,N)}{Diff^+(M)}$$
we get:
\begin{Proposition}
$B_+(M,N)$ is a 2-cover of $B(M,N).$
\end{Proposition}

Now, taking  basepoints $x_0 \in M$ and $y_0 \in N,$ we define the principal bundle of based embeddings

\begin{Proposition}
Let $$Emb_b(M,N) = \{f \in Emb(M,N) | f(x_0) = y_0 \}.$$
Let $$Diff_b(M) = \{ g \in Diff(M) | g(x_0) = x_0 \}.$$
Let $$ Diff_{b,+}(M) = Diff_b(M) \cap Diff_+(M).$$
Let $$B_b(M,N) = Emb_b(M,N) / Diff_b(M,N)$$
and $$B_{b,+}(M,N) = Emb_{b}(M,N) / Diff_{b,+}(M,N).$$
Then $Emb_b(M,N)$ is a principal bundle with base $B_b(M,N)$ (resp. $B_{b,+}(M,N)$) and with structure group $Diff_b(M)$ (resp. $Diff_{b,+}(M)$) 
\end{Proposition}

\noindent
\textbf{Proof.} It follows from the fact that $Emb_b(M,N) = ev_{x_0}^{-1}(y_0)$ in $Emb(M,N),$ and $Diff_b(M) = ev_{x_0}^{-1}(x_0)$ in $Diff(M).$ \qed

\section{Chern-Weil forms on principal bundle of embeddings and homotopy invariants}
\subsection{Chern forms in infinite dimensional setting}

Let $P$ be a
principal bundle, of basis $M$ and with structure group $G$. Let
$\mathfrak g$ be the Lie algebra of $G$. Recall that $G$ acts on
$P$, and also on $P \times \mathfrak g$ by the action $((p,v), g)
\in (P \times \mathfrak g)\times G \mapsto (p.g, Ad_{g^{-1}}(v))
\in (P \times \mathfrak g)$. Let $AdP = P \times_{Ad g} = (P \times
\mathfrak g )/ G $ be the adjoint bundle of $P$, of basis $M$ and
of typical fiber $\mathfrak g$, and let $Ad^kP = (Ad P)^{\times
k}$ be the product bundle, of basis $M$ and of typical fiber
$\mathfrak g^{\times k}$.

\begin{Definition} Let $k$ in $\N^*$. We define
$\pol^{k}(P)$, the set of smooth maps $Ad^kP \rightarrow \C$ that are $k$-linear
and symmetric on each fiber, equivalently as the set of smooth maps 
$P \times \mathfrak{g}^k
 \rightarrow \C$ that are $k$-linear symmetric in the second variable 
and $G$-invariants with respect to the natural coadjoint action 
of $G$ on $\mathfrak{g}^k .$

Let $\pol(P) = \bigoplus_{k \in \N*}\pol(P)$.
\end{Definition}

Let $\Con(P)$ be the set of connections on $P$. For any $\theta \in \Con(P)$,
we denote by $F(\theta)$ its curvature and $\nabla^\theta$ 
(or $\nabla$ when it carries no ambiguity) its covariant derivation. 
Given an algebra $A$,
In this section, we study the maps, for $k \in \N^*$,

\begin{eqnarray} \label{fomega}
 Ch \quad : \quad \Con(P) \times \pol^k(P) & \rightarrow & \Omega^{2k}(M,\C)  \\
 (\theta, f) & \mapsto& Alt (f(F(\theta),...,F(\theta))) \end{eqnarray}
where $Alt$ denotes the skew-symmetric part of the form.
Notice that, in the case of the finite dimensional matrix groups
$Gl_n$ with Lie algebra $\mathfrak{gl}_n$, the set $\pol(P)$ 
is generated by the polynomials $A \in \mathfrak{gl}_n \mapsto \hbox{tr}(A^k),$ for $
k \in {0,...,n}$. This leads to classical definition of Chern forms. 
However, in the case of infinite dimensional structure groups, most situations are 
still unknown and we do not know how to define a set of generators for $\pol(P).$

\begin{Lemma} \label{bracket}
Let $f \in \pol^k(P).$ Then \begin{eqnarray*} \label{1} f([a_1,v], a_2, ..., a_k ) + f(a_1, [a_2,v],...,a_k) + \\
... && \\
+ f(a_1, a_2,...,[a_k,v])&=& 0. \end{eqnarray*}
\end{Lemma}

\vskip 12pt
\noindent
\textbf{Proof.}
Let us notice first that $f$ is symmetric. 
Let $v \in \mathfrak{g},$ and $c_t$ a path in $G$ such that 
$\lbrace {d \over dt}c_t \rbrace_{t=0}=v.$
Let $a_1,...,a_k \in \mathfrak{g}^k$. 
\begin{eqnarray*}
\lbrace {d \over dt}\lbrace f(ad_{c_t^{-1}}a_1,...,ad_{c_t^{-1}}a_k)\rbrace_{t=0}
& = &  f([a_1,v], a_2, ..., a_k ) + f(a_1, [a_2,v],...,a_k) + ... 
+ f(a_1, a_2,...,[a_k,v])\end{eqnarray*}
Since $f$ in $G$-invariant, we get 
\begin{eqnarray*} \label{1'} f([a_1,v], a_2, ..., a_k ) + f(a_1, [a_2,v],...,a_k) + \\
... && \\
+ f(a_1, a_2,...,[a_k,v])&=& 0. \end{eqnarray*}

\begin{Lemma} \label{local}
Let $f\in \pol^k(P)$ such that $f$, as 
a smooth map $P\times \mathfrak{g}^k \rightarrow \C,$ satifies $d^Mf=
0$ on a system of local trivializations of $P.$ 
Then, the map 
$$Ch^f : \theta \in \Con(P) \mapsto Ch^f(\theta)= Ch(\theta,f)\in \Omega^*(P,\C)$$
takes values into closed forms on $P$. Moreover, 

(i) it is vanishing on vertical 
vectors and defines a closed form on $M$.

(ii) the cohomology class of this form does not depend on the choice of the
chosen connexion $\theta$ on $P$. 
\end{Lemma}

\noindent
\textbf{Proof.}
The proof runs as in the finite dimensional case, see e.g.\cite{KN}
First, it is vanishing on vertical vectors and $G$-invariant because 
the curvature of a connexion
vanishes on vertical forms and is $G$-covariant for the coadjoint action.
Let us now fix $f \in \pol^k(P)$. 
We compute       
$df(F(\theta),...,F(\theta)).$ We notice first that it vanishes on vertical vectors trivially. 
Let us fix $Y_1^h,...,Y_{2k}^h,X^h$ $2k+1$ horizontal vectors on $P$ at $p \in P$.
On a local trivialization of $P$ around $p$, these vectors read as 
\begin{eqnarray*}
Y_1^h & = & Y_1 - \tilde \theta(Y_1) \\
 & (...) & \\
Y_{2k}^h & = & Y_{2k} - \tilde \theta(Y_{2k}) \\
X^h & = & X - \tilde \theta(X) \end{eqnarray*}
where $\tilde \theta$ stands here for the expression of $\theta$ in the local 
trivilization, and $Y_1,...,Y_{2k},X$ $2k+1$ tangent vectors on $M$ at $\pi(p) \in M.$
We extend these vector fields on a neighborhood of $p$

- by the action of $G$ in the vertical directions

- setting the vectors fields 
constant on $U \times p$, where $U$ is a local chart on $M$ 
around $\pi(p)$.

Then, we have
$$ f(F(\theta),...,F(\theta))(Y_1^h,...,Y_{2k}^h) = 
f(F(\theta),...,F(\theta))(Y_1,...,Y_{2k})$$   
since $F(\theta)$ is vanishing on vertical vectors.

Then, on a local trivialization 
with the notations defined before
(the sign $Alt$ is omitted for easier reading),
and writing $d^M$ for the differential of forms on any open subset of $M$,       

\begin{eqnarray*}
d^Mf(F(\tilde \theta),...,F(\tilde\theta)) & = & \sum_{i = 1}^k f(d^MF(\tilde\theta),F(\tilde\theta)
,...,F(\tilde\theta)) + f(F(\tilde\theta), d^MF(\tilde\theta),...,F(\tilde\theta)) + \\
&& ... 
+ f(F(\tilde\theta), F(\tilde\theta),..., d^MF(\tilde\theta))\end{eqnarray*}

and then, using Lemma \ref{bracket},
\begin{eqnarray*}
\nabla^\theta f(F(\tilde \theta),...,F(\tilde\theta)) & = & \sum_{i = 1}^k f(\nabla^\theta F(\tilde\theta),F(\tilde\theta)
,...,F(\tilde\theta)) + f(F(\tilde\theta), \nabla^\theta F(\tilde\theta),..., F(\tilde\theta)) + \\
&& ... 
+ f(F(\tilde\theta), F(\tilde\theta),..., \nabla^\theta F(\tilde\theta))
\end{eqnarray*}
Then, by Bianchi identity, 
we get that 
\begin{eqnarray*}
d^MCh(f,\theta) & = & \nabla^\theta Ch(f,\theta)\\
& = & 0 \end{eqnarray*}
This proves (i)
Then, following e.g \cite{KN}, if $\theta$ and $\theta'$ are connections, fix 
$\mu = \theta' - \theta$ and $\theta_t = \theta + t\nu$ for $t \in [0;1].$
We have \begin{eqnarray*}
{d F(\theta_t) \over dt} &=& \nabla^{\theta^t} \mu\end{eqnarray*}
Moreover, 
$\mu$ is $G$-invariant and vanishes on vertical vectors. Thus, 
\begin{eqnarray*}
{d Ch(f,\theta_t) \over dt}& = &k f(F(\theta_t),...,F(\theta_t),\nabla^{\theta_t}\mu)\\
& = & k d^M (f(F(\theta_t),...,F(\theta_t),\mu)). \end{eqnarray*} 
Integrating in the $t$-variable, we get
$$Ch(f,\theta_0) - Ch(f,\theta_1) = 
-kd^M \int_0^1 f(F(\theta_t),...,F(\theta_t),\mu) dt.$$
Even if these computations are local, the two sides are global objects 
and do not depend on the chosen trivialization, which ends the proof.
\qed

\vskip 12pt
\noindent
\textbf{Important remark.}
The condition $d^Mf=0$ is a \textbf{local} condition, checked in an 
(adequate) system of trivializations of the principal bundle, 
because it has to be checked on the vector bundle $Ad(P)^{\times k}.$ 
This is in 
particular the case when we can find a $0$-curvature 
connection $\theta$ on $P$ such that
$$[\nabla^\theta, f] = 0$$
In that case, since the structure group $G$ is regular, 
we can find a system of local trivializations of $P$ defined by $\theta$
and such that,on any local trivialization, $\nabla^\theta = d^M$  
(see e.g. \cite{KM}, \cite{Ma2004} for the technical tools that are
necessary for this).
\vskip 12pt
\noindent
This technical remark can appear rather unsatisfactory 
first because it restricts the 
ability of application of the previous lemma, secondly because we need 
have a local (and rather unelegant) condition. 
This is why we give the following theorem, from Lemma \ref{local}. 
\begin{Theorem} \label{chern-weil}
Let $f\in \pol(P)$ for which there exists $\theta \in \mathcal{C}(P)$ such that $[\nabla^\theta,f]=0.$ 
We shall note this set of polynomials by $\pol_{reg}(P).$
Then, the map 
$$Ch^f : \theta \in \Con(P) \mapsto Ch^f(\theta)= Ch(\theta,f)\in \Omega^*(P,\C)$$
takes values into closed forms on $P$. Moreover, 

(i) it is vanishing on vertical 
vectors and defines a closed form on $M$.

(ii) the cohomology class of this form does not depend on the choice of the
chosen connexion $\theta$ on $P$. 

Moreover,
$\forall (\theta,f) \in\mathcal{C}(P)\times \pol_{reg}(P), [\nabla^\theta,f] = 0.$
\end{Theorem}

\noindent
\textbf{Proof.}
Let $f \in \pol_{reg}(P)$ and let $\theta \in \mathcal{C}(P)$ such that $[\nabla^\theta , f]=0.$ 
Let $\theta' \in \mathcal{c}(P)$ and 
let $\nu = \theta' - \theta \in \Omega^1(M,\mathfrak{g}).$
Let $(\alpha_1,...,\alpha_k)\in (\Omega^2(M,\mathfrak{g}))^k.$

\begin{eqnarray*}
[\nabla^{\theta'},f](\alpha_1,...\alpha_k) & = &
[\nabla^{\theta},f](\alpha_1,...\alpha_k) + f([\alpha_1,\nu],...,\alpha_n)+ \\
&& ... + f(\alpha_1,...,[\alpha_n,\nu])\\
& = & f([\alpha_1,\nu],...,\alpha_n)+... + f(\alpha_1,...,[\alpha_n,\nu]) \\
& = & 0.
\end{eqnarray*} 
Then, 
$\forall (\theta,f) \in\mathcal{C}(P)\times \pol_{reg}(P), [\nabla^\theta,f] = 0.$
By the way, $\forall \theta' \in \mathcal{C}(P),$ 
$$d^M f(\alpha_1,...,\alpha_k) = f(\nabla^{\theta'} \alpha_1,...,\alpha_k)+
... + f(\alpha_1,...,\nabla^{\theta'} \alpha_k).$$
Applying this to $\alpha_1=...=\alpha_k = F(\theta'),$
we get $$ dCh(f,\theta') = f(\nabla^{\theta'} F(\theta'),...,F(\theta'))+
... + f(F(\theta'),...,\nabla^{\theta'} F(\theta')) = 0$$
by Bianchi identity. Thus $Ch(f,\theta')$ is closed.
Then, mimicking the end of the proof of Lemma \ref{local}, we get that the
difference $Ch(f,\theta) - Ch(f,\theta')$ is an exact form, which ends the proof.

\begin{Proposition}
Let $\phi : \mathfrak{g}^k \rightarrow \C$ be a $k-$linear, symmetric, $Ad-$invariant form. Let $f:P\times \mathfrak{g}^k\rightarrow \C$ be the map induced by $\phi$ by the formula:
$f(x,g)= \phi(g).$ 
Then $f \in \pol_{reg}.$ 
\end{Proposition}

\vskip 12pt
\noindent
\textbf{Proof.}
Obsiously, $f\in \pol .$ 
Let $\varphi: U\times G \rightarrow P$ and $ \varphi ': U \times G \rightarrow P$ be a local trivialisations of $P$, where $U$ is an open subset of $M.$
Then there exists a smooth map $g : U \rightarrow G$ such that  $\varphi '(x, e_G) = \varphi(x,e_G).g(x) .$ Then we remark that $\varphi^* f = \varphi '^* f $
is a constant map on horizontal slices since $\phi$ is Ad-invariant. Moreover,
since  $\varphi^* f$ in a constant (polynomial-valued) map on  $\varphi(x,e_G)$
we get that $[\nabla^\theta , f] = 0$ for the (flat) connection $\theta$ such that $T\varphi(x,e_G)$ spans the horizontal bundle over $U$. \qed

\subsection{Application to $Emb(M,N)$}
Mimicking the approach of \cite{Ma2006}, the cohomology classes of Chern-Weil forms should give rise to homotopy invariants.  
Applying Theorem \ref{chern-weil}, we get:

\begin{Theorem} \label{cantrace}
The Chern-Weil forms $Ch^f$ is a $H^*(B(M,N))-$valued invariant of the homotopy class of an embedding, $\forall k \in \N^*.$ 
\end{Theorem}

When $M=S^1$, $Emb(S^1,N)$ is the space of (parametrized) smooth knots on $N$, and $B(S^1,N)$ is the
space of non parametrized knots. Its connected components are the homotopy classes
of the knots, through classical results of differential topology, see e.g. \cite{Hir}.
We now apply the material of the previous section to manifolds of embeddings. For this, we can define 
invariant polynomials of the type of those obtained in \cite{Ma2006} (for mapping spaces) by a field of linear functionnal $\lambda$ with ``good properties'' that ensures that 
$$ A \mapsto \lambda(A^k) \in \pol^k_{reg}.$$ 
This approach is a straightforward generalization of the description of Chern-Weil forms on finite dimensional principal bundles where polynomials are generated by functionnals of the type $A \mapsto \tr(A^k)$ ($\tr$ is the classical trace) but as we guess that we  can consider other classes of polynomials for spaces of embeddings. In 0this paper, let us describe how to replace the classical trace of matrices $\tr$ by a renormalized trace $\tr^Q.$
In the most general case, it is not so easy to define a family of weights $f \in Emb(M,N) \mapsto Q_f$ which satisfy the good properties. Indeed, we have two examples of constructions which match the necessary assumptions for $\pol_{reg}$ when $M= S^1$, and the first one is derived from  the following example:

\subsubsection*{Knot invariant through Kontsevich and Vishik trace}

The Kontsevich and Vishik trace is a renormalized trace for which $tr^Q([A,B]) = 0$ for each differential operator $A,B$ and does not depend 
on the weight chosen in the odd class. For example, one can choose $Q = Id+  {\nabla } ^* {\nabla  }$, where $\nabla$ is a connection induced on $\mathcal{N}_f$ by the Riemannian metric, as described in \cite{Ma2006}. It is an order 2 injective elliptic differential operator (in the odd class), and the coadjoint action of $Aut(\mathcal{N}_f)$ will give rise to another order 2 injective elliptic differential operator \cite{Gil}. When $Q = Id+  {\nabla } ^* {\nabla  },$ this only changes $\nabla$ into another connection on $E.$ Thus, setting 
$$\phi(A,...,A) = tr^Q(A^k), $$
we have $$f \in \pol_{reg}.$$
Let us now consider a connected component of $B(M,N),$ i.e. a homotopy class of an embedding among the space of embeddings.
We apply now the construction to $M=S^1.$
The polynomial 
$$ \phi :A \mapsto tr^Q(A^k)$$ is $Diff(S^1)-$invariant, and gives rise to an invariant of non oriented knots, i.e. a Chern form on the base manifold $$B(S^1,N) = Emb(S^1,N) / Diff(S^1)$$
by theorem \ref{cantrace}.
This approach can be extended to invanraint of embeddings, replacing $S^1$ by another odd-dimensional manifold.

\end{document}